\documentclass{article}                                     
                                                                             
\input{amssym.def}                                                           
\input{amssym.tex}                                                           
                                                                             
\begin{document}                                                             
\title{Remark on the rank of elliptic curves}

\author{Igor  ~Nikolaev
\footnote{Partially supported 
by NSERC.}}


\date{}
 \maketitle


\newtheorem{thm}{Theorem}
\newtheorem{lem}{Lemma}
\newtheorem{dfn}{Definition}
\newtheorem{rmk}{Remark}
\newtheorem{cor}{Corollary}
\newtheorem{cnj}{Conjecture}
\newtheorem{exm}{Example}


\newcommand{\N}{{\Bbb N}}
\newcommand{\F}{{\cal F}}
\newcommand{\R}{{\Bbb R}}
\newcommand{\Z}{{\Bbb Z}}
\newcommand{\C}{{\Bbb C}}

\begin{abstract}
A covariant functor on the elliptic curves with complex multiplication
is constructed. The functor takes values in the noncommutative
tori with real multiplication.  A conjecture on  the rank of an elliptic curve 
is formulated.

\vspace{7mm}


\vspace{5mm}
{\it AMS (MOS) Subj. Class.: 11G15; 46L85}
\end{abstract}

\section*{Introduction}
{\bf A.} Let $0<\theta< 1$ be an irrational number given by the regular
continued fraction 
$$
\theta=a_0+{1\over\displaystyle a_1+
{\strut 1\over\displaystyle a_2\displaystyle +\dots}}
=[a_0, a_1, a_2,\dots].
$$
Consider an $AF$-algebra, ${\Bbb A}_{\theta}$,  defined by the Bratteli 
diagram:

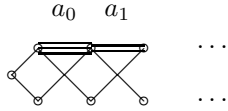
\begin{figure}[here]
\begin{picture}(300,60)(0,0)
\put(110,30){\circle{3}}
\put(120,20){\circle{3}}
\put(140,20){\circle{3}}
\put(160,20){\circle{3}}
\put(120,40){\circle{3}}
\put(140,40){\circle{3}}
\put(160,40){\circle{3}}

\put(110,30){\line(1,1){10}}
\put(110,30){\line(1,-1){10}}
\put(120,42){\line(1,0){20}}
\put(120,40){\line(1,0){20}}
\put(120,38){\line(1,0){20}}
\put(120,40){\line(1,-1){20}}
\put(120,20){\line(1,1){20}}
\put(140,41){\line(1,0){20}}
\put(140,39){\line(1,0){20}}
\put(140,40){\line(1,-1){20}}
\put(140,20){\line(1,1){20}}

\put(180,20){$\dots$}
\put(180,40){$\dots$}

\put(125,52){$a_0$}
\put(145,52){$a_1$}

\end{picture}

\caption{The $AF$-algebra  ${\Bbb A}_{\theta}$.}
\end{figure}

\noindent
where $a_i$ indicate the multiplicity  of the edges of the graph.
(For a definition of the $AF$-algebras and their Bratteli diagrams,
we refer the reader to \cite{E}, or \S 1.2.)
For the simplicity, we shall say that ${\Bbb A}_{\theta}$ is a noncommutative torus.
Note that the classical definition of a noncommutative torus is slightly  different but equivalent from
the standpoint of the $K$-theory \cite{EfSh1},  \cite{PiVo1}, \cite{Rie1}. 
The ${\Bbb A}_{\theta}$ is said to have real multiplication, if $\theta$
is a quadratic irrationality. 
Recall that the noncommutative tori ${\Bbb A}_{\theta}, {\Bbb A}_{\theta'}$ are 
stably isomorphic whenever ${\Bbb A}_{\theta}\otimes {\cal K}
\cong {\Bbb A}_{\theta'}\otimes {\cal  K}$, where ${\cal K}$ is the $C^*$-algebra
of the compact operators. It is known that ${\Bbb A}_{\theta}, {\Bbb A}_{\theta'}$ are stably
isomorphic if and only if $\theta'\equiv \theta~mod~GL(2, {\Bbb Z})$,
i.e. $\theta'= (a\theta +b)~/~ (c\theta +d)$,  where $a,b,c,d\in {\Bbb Z}$ and $ad-bc=\pm 1$.

\medskip
{\bf B.}
Let  $\Lambda=\omega_1{\Bbb Z}+\omega_2{\Bbb Z}$ be a lattice in the complex plane ${\Bbb C}$.
Recall that $\Lambda$  defines an elliptic curve $E({\Bbb C}): y^2=4x^3-g_2x-g_3$
via the complex analytic map ${\Bbb C}/\Lambda\to E({\Bbb C})$ given by
the formula $z\mapsto (\wp (z,\Lambda), \wp'(z,\Lambda))$, where  
$g_2=60\sum_{\omega\in\Lambda^{\times}}\omega^{-4}$,
$g_3=140\sum_{\omega\in\Lambda^{\times}}\omega^{-6}$,
$\Lambda^{\times}=\Lambda-\{0\}$ and 
$$
 \wp (z,\Lambda)={1\over z^2}+\sum_{\omega\in\Lambda^{\times}} \left({1\over (z-\omega)^2}-{1\over \omega^2}\right)
$$
 is the Weierstrass $\wp$ function. We identify the elliptic curves $E({\Bbb C})$ 
with the complex tori ${\Bbb C}/\Lambda$. If $\tau=\omega_2/\omega_1$,
then $E_{\tau}({\Bbb C}), E_{\tau'}({\Bbb C})$ are isomorphic whenever  
$\tau'\equiv \tau ~mod~GL(2, {\Bbb Z})$.
The endomorphism 
ring $End~({\Bbb C}/\Lambda)$ is isomorphic either to ${\Bbb Z}$ or to
an order in the imaginary quadratic number field $k$ \cite{S}. In the second case,
we say that the elliptic curve has a complex multiplication and denote such a curve
by $E_{CM}$.

\medskip
{\bf C.}
Consider the cubic $E_{\lambda}: y^2=x(x-1)(x-\lambda)$, $\lambda\in {\Bbb C}-\{0,1\}$.
The $j$-invariant of $E_{\lambda}$ is given by the formula
$j(E_{\lambda})=2^6 (\lambda^2-\lambda+1)^3\lambda^{-2}(\lambda-1)^{-2}$.
To find $\lambda$ corresponding to an elliptic curve with the complex multiplication,
one has to solve the polynomial equation $j(E_{CM})=j(E_{\lambda})$ with respect
to $\lambda$. Since $j(E_{CM})$ is an algebraic integer, $\lambda_{CM}\in K$,
where $K$ is an algebraic extension (of the degree at most six) of the field
${\Bbb Q}(j(E_{CM}))$. Thus, each $E_{CM}$ is isomorphic to the cubic
$y^2=x(x-1)(x-\lambda_{CM})$ defined over the field $K$. 
The Mordell-Weil theorem  says that the set of the $K$-rational
points of $E_{CM}$ is a finitely generated abelian group, 
whose rank we shall  denote  by $rk~(E_{CM})$.

\medskip
{\bf D.} Let ${\cal E}$ be a category whose objects are elliptic curves
and the arrows are isomorphisms of the elliptic curves. Likewise, let ${\cal A}$
be a category whose objects are noncommutative tori and the arrows are
stable isomorphisms of the noncommutative tori. Our main goals can be 
expressed as follows. 

\medskip\noindent
{\bf Objectives.} 
{\it (i) to construct a functor (if any) $F: {\cal E}\to {\cal A}$,  which
maps isomorphic elliptic curves to the stably isomorphic noncommutative
tori;   (ii) to study the range of $F$ on the elliptic curves with complex multiplication and
(iii) to interpret  the invariants of the stable isomorphism
classes of the noncommutative tori in terms of the arithmetic invariants of the elliptic 
curves.}

\medskip\noindent
In the course of this note, we were able to obtain an  answer to (i) and (ii), 
while (iii) generates a conjecture. Namely, a covariant 
non-injective functor $F: {\cal E}\to {\cal A}$, which maps isomorphic
elliptic curves to the stably isomorphic noncommutative tori,  is
constructed (lemma \ref{lm1}). It is proved that $F$ sends the elliptic 
curves with complex multiplication to the noncommutative tori
with real multiplication (theorem \ref{thm1}). Finally, a conjecture
on the rank of an elliptic curve with the complex multiplication is
formulated (\S 3). The functor $F$ has been studied by Kontsevich \cite{Kon1} (e.g. \S 1.39),
Manin \cite{Man1},  Polishchuk \cite{Pol1}-\cite{Pol4}, 
Polishchuk-Schwarz \cite{PoSch1}, Soibelman \cite{Soi1} and \cite{Soi2}, 
Soibelman-Vologodsky \cite{SoVo1}, Taylor \cite{Tay1} and \cite{Tay2} 
{\it et al.}  Our terminology  is freely and gratefully  borrowed from the 
above works.

\medskip
{\bf E.}
The existence and properties of $F$ are part of a Hodge theory for the measured
foliations on a closed surface. Such  a theory has been developed
by Hubbard and Masur \cite{HuMa1}, who were inspired by the works of Thurston \cite{Thu1}.
We shall give in \S 1 a brief account of the Hubbard-Masur-Thurston theory  
and the explicit formulas for the functor $F$. At the heart of the construction 
is a diagram:
$$   
F: {\cal E}
\buildrel\rm h\over
\longrightarrow
{\Bbb R}^2
\buildrel\rm\pi \over
\longrightarrow
{\Bbb R}P^1\cong {\cal A},
$$
where $h$ is a bijection and $\pi$ is a projection map. 
For the sake of brevity,  let $Isom~(E)=\{E'\in {\cal E}~|~E'\cong E\}$
be the isomorphism class of an elliptic curve $E$, 
$h(Isom~(E))=\mu_E({\Bbb Z}+\theta_E{\Bbb Z}):={\goth m}_E\subset {\Bbb R}$ be a ${\Bbb Z}$-module
and $F(E)={\Bbb A}_{\theta_E}$. A summary of our
results can be formulated as follows. 
\begin{lem}\label{lm1}
Let $\varphi: E\to E'$ be an isogeny of the elliptic curves. 
Then $\theta'\equiv\theta~mod~M_2({\Bbb Z})$, where $M_2({\Bbb Z})$
is an integer matrix of the rank 2. In particular, $F$ maps the isomorphic 
elliptic curves  to the stably isomorphic noncommutative tori. 
\end{lem}
\begin{thm}\label{thm1}
Let $E\in Isom~(E_{CM})$. Then there exists an $h$, such that:

\medskip
(i) ${\goth m}_E$ is a full module in the  real quadratic number field;

\smallskip
(ii) ${\goth m}_E$ is an invariant of the class $Isom~(E_{CM})$.

\medskip\noindent
In particular, $\theta_E$ is a quadratic irrationality.
\end{thm}
The structure of the note  is as follows. In section 1,  we 
introduce the notation and some preliminary facts. 
The lemma \ref{lm1} and theorem \ref{thm1} are 
proved in the section 2. In section 3,  a conjecture on the  rank
of an elliptic curve is formulated.

\section{Preliminaries}
This section contains a  summary of measured foliations, 
$AF$-algebras and the functor $F$. 
The reader is encouraged to consult \cite{E} (operator algebras)  and  
\cite{HuMa1} (measured foliations \& Teichm\"uller space) for a 
systematic account.

\subsection{Measured foliations and $T(g)$}
{\bf A.}  A measured foliation, ${\cal F}$, on a surface $X$
is a  partition of $X$ into the singular points $x_1,\dots,x_n$ of
order $k_1,\dots, k_n$ and the regular leaves (1-dimensional submanifolds). 
On each  open cover $U_i$ of $X-\{x_1,\dots,x_n\}$ there exists a non-vanishing
real-valued closed 1-form $\phi_i$  such that 

\medskip
(i)  $\phi_i=\pm \phi_j$ on $U_i\cap U_j$;

\smallskip
(ii) at each $x_i$ there exists a local chart $(u,v):V\to {\Bbb R}^2$
such that for $z=u+iv$, it holds $\phi_i=Im~(z^{k_i\over 2})$ on
$V\cap U_i$ for some branch of $z^{k_i\over 2}$. 

\medskip\noindent
The pair $(U_i,\phi_i)$ is called an atlas for the measured foliation ${\cal F}$.
Finally, a measure $\mu$ is assigned to each segment $(t_0,t)\in U_i$, which is  transverse to
the leaves of ${\cal F}$, via the integral $\mu(t_0,t)=\int_{t_0}^t\phi_i$. The 
measure is invariant along the leaves of the foliation ${\cal F}$, hence the name.

\smallskip
{\bf B.} Let $S$ be a Riemann surface, and $q\in H^0(S,\Omega^{\otimes 2})$ a holomorphic
quadratic differential on $S$. The lines $Re~q=0$ and $Im~q=0$ define a pair
of measured foliations on $R$, which are transversal to each other outside the set of 
singular points. The set of singular points is common to the both foliations and coincides
with the zeroes of $q$. The above measured foliations are said to represent the  vertical and horizontal 
 trajectory structure  of $q$, respectively.

\smallskip
{\bf C.}  Let $T(g)$ be the Teichm\"uller space of the topological surface $X$ of genus $g$,
i.e. the space of complex structures on $X$. 
Consider the vector bundle $p: Q\to T(g)$ over $T(g)$ whose fiber above a point 
$S\in T(g)$ is the vector space $H^0(S,\Omega^{\otimes 2})$.   
Given non-zero $q\in Q$ above $S$, one can consider the horizontal measured foliation
${\cal F}_q\in \Phi_X$ of the quadratic differential $q$, where $\Phi_X$ is  the space of 
(equivalence classes of) 
measured foliations on $X$. If $\{0\}$ is the zero section of $Q$,
the above construction defines a map $Q-\{0\}\longrightarrow \Phi_X$. 
For any ${\cal F}\in\Phi_X$, let $E_{\cal F}\subset Q-\{0\}$ be the fiber
above ${\cal F}$. In other words, $E_{\cal F}$ is a subspace of the holomorphic 
quadratic differentials, whose horizontal trajectory structure coincides with the 
measured foliation ${\cal F}$. 

\bigskip\noindent
{\bf Theorem (\cite{HuMa1})}
{\it The restriction $E_{\cal F}\longrightarrow  T(g)$ of $p$ to $E_{\cal F}$ is a
homeomorphism.}

\medskip
{\bf D.} 
Let  $\Phi_X$ be  the space  of measured
foliations on the topological surface $X$. Following Douady and Hubbard
\cite{DoHu1}, we shall consider a coordinate system on $\Phi_X$,
suitable for the construction of the functor $F$. 
For clarity, let us make a generic assumption that $q\in H^0(S,\Omega^{\otimes 2})$
is a holomorphic quadratic differential with the simple zeroes only. 
We wish to construct a Riemann surface of $\sqrt{q}$, which is a double cover
of $S$ with the ramification over the zeroes of $q$. Such a surface, denoted by
$\widetilde S$, is unique and has an advantage of carrying a holomorphic
differential $\omega$, such that $\omega^2=q$. Denote by 
$\pi:\widetilde S\to S$ a covering projection. The vector space
$H^0(\widetilde S,\Omega)$ splits into the direct sum
$H^0_{even}(\widetilde S,\Omega)\oplus H^0_{odd}(\widetilde S,\Omega)$
in view of  the involution $\pi^{-1}$ of $\widetilde S$, and
the vector space $H^0(S,\Omega^{\otimes 2})\cong H^0_{odd}(\widetilde S,\Omega)$.
Let $H_1^{odd}(\widetilde S)$ be an odd part of the homology of $\widetilde S$
relatively  the zeroes of $q$.   Consider a pairing
$H_1^{odd}(\widetilde S)\times H^0(S, \Omega^{\otimes 2})\to {\Bbb C}$,
defined by the integration  $(\gamma, q)\mapsto \int_{\gamma}\omega$. 
Take the associated map
$\psi_q: H^0(S,\Omega^{\otimes 2})\to Hom~(H_1^{odd}(\widetilde S); {\Bbb C})$
and let $h_q= Re~\psi_q$. 

\bigskip\noindent
{\bf Theorem (\cite{DoHu1})}
{\it The map
$h_q: H^0(S, \Omega^{\otimes 2})\longrightarrow Hom~(H_1^{odd}(\widetilde S); {\Bbb R})$
is an ${\Bbb R}$-isomorphism.} 

\medskip\noindent
Since  each  ${\cal F}\in \Phi_X$ is the  vertical foliation 
$Re~q=0$ for a $q\in H^0(S, \Omega^{\otimes 2})$, the theorem
implies that $\Phi_X\cong Hom~(H_1^{odd}(\widetilde S); {\Bbb R})$.
By the  formulas for the relative homology: 
$$
H_1^{odd}(\widetilde S)\cong {\Bbb Z}^n, ~\hbox{where}
~n=\cases{6g-6, & if $g\ge 2$\cr 2, & if $g=1$.}     
$$
Thus, if $\{\gamma_1,\dots,\gamma_n\}$ is a basis in 
$H_1^{odd}(\widetilde S)$, the reals $\lambda_i=\int_{\gamma_i} Re~\omega$
are natural coordinates in the space $\Phi_X$ \cite{DoHu1}.

\subsection{AF-algebras}
{\bf A.}
The $C^*$-algebra is an algebra $A$ over ${\Bbb C}$ with a norm
$a\mapsto ||a||$ and an involution $a\mapsto a^*$ such that
it is complete with respect to the norm and $||ab||\le ||a||~ ||b||$
and $||a^*a||=||a^2||$ for all $a,b\in A$.
If $A$ is commutative, then the Gelfand  theorem says that $A$ is isomorphic
to the $C^*$-algebra $C_0(X)$ of continuous complex-valued
functions on a locally compact Hausdorff space $X$.
For otherwise, $A$ represents a noncommutative  topological
space $X$.

\medskip
{\bf B.} Let $A$ be a $C^*$-algebra deemed as a noncommutative
topological space. One can ask when two such topological spaces
$A,A'$ are homeomorphic? To answer the question, let us recall
the topological $K$-theory.  If $X$ is a (commutative) topological
space, denote by $V_{\Bbb C}(X)$ an abelian monoid consisting
of the isomorphism classes of the complex vector bundles over $X$
endowed with the Whitney sum. The abelian monoid $V_{\Bbb C}(X)$
can be made to an abelian group, $K(X)$, using the Grothendieck
completion. The covariant functor $F: X\to K(X)$ is known to
map the homeomorphic topological spaces $X,X'$ to the isomorphic
abelian groups $K(X), K(X')$.  Let now $A,A'$ be the $C^*$-algebras. If one wishes to
define a homeomorphism between the noncommutative topological spaces $A$ and $A'$, 
it will suffice to define an isomorphism between the abelian monoids $V_{\Bbb C}(A)$
and $V_{\Bbb C}(A')$ as suggested by the topological $K$-theory. 
The r\^ole of the complex vector bundle of degree $n$ over the 
$C^*$-algebra $A$ is played by a $C^*$-algebra $M_n(A)=A\otimes M_n$,
i.e. the matrix algebra with the entries in $A$.  The abelian monoid
$V_{\Bbb C}(A)=\cup_{n=1}^{\infty} M_n(A)$ replaces the monoid 
$V_{\Bbb C}(X)$ of the topological $K$-theory. Therefore, 
the noncommutative topological spaces $A,A'$ are homeomorphic,
if $V_{\Bbb C}(A)\cong V_{\Bbb C}(A')$ are isomorphic abelian
monoids. The latter equivalence is called a {\it stable isomorphism}
of the $C^*$-algebras $A$ and $A'$ and is formally written as 
$A\otimes {\cal K}\cong A'\otimes {\cal K}$, where 
${\cal K}=\cup_{n=1}^{\infty}M_n$ is the $C^*$-algebra of compact
operators.  Roughly speaking, the stable isomorphism between
the $C^*$-algebras $A$ and $A'$ means that $A$ and $A'$ are
homeomorphic as the noncommutative topological spaces.

\medskip
{\bf C.}
Let $A$ be a unital $C^*$-algebra and $V(A)$
be the union (over $n$) of projections in the $n\times n$
matrix $C^*$-algebra with entries in $A$.
Projections $p,q\in V(A)$ are equivalent if there exists a partial
isometry $u$ such that $p=u^*u$ and $q=uu^*$. The equivalence
class of projection $p$ is denoted by $[p]$.
The equivalence classes of orthogonal projections can be made to
a semigroup by putting $[p]+[q]=[p+q]$. The Grothendieck
completion of this semigroup to an abelian group is called
a  $K_0$-group of algebra $A$.
Functor $A\to K_0(A)$ maps a category of unital
$C^*$-algebras into the category of abelian groups so that
projections in algebra $A$ correspond to a positive
cone  $K_0^+\subset K_0(A)$ and the unit element $1\in A$
corresponds to an order unit $u\in K_0(A)$.
The ordered abelian group $(K_0,K_0^+,u)$ with an order
unit  is called a {\it dimension group}.

\medskip
{\bf D.}
An {\it $AF$-algebra}  (approximately finite $C^*$-algebra) is defined to
be the  norm closure of an ascending sequence of the finite dimensional
$C^*$-algebras $M_n$'s, where  $M_n$ is the $C^*$-algebra of the $n\times n$ matrices
with the entries in ${\Bbb C}$. Here the index $n=(n_1,\dots,n_k)$ represents
a semi-simple matrix algebra $M_n=M_{n_1}\oplus\dots\oplus M_{n_k}$.
The ascending sequence mentioned above  can be written as 
$$M_1\buildrel\rm\varphi_1\over\longrightarrow M_2
   \buildrel\rm\varphi_2\over\longrightarrow\dots,
$$
where $M_i$ are the finite dimensional $C^*$-algebras and
$\varphi_i$ the homomorphisms between such algebras.  The set-theoretic limit
$A=\lim M_n$ has a natural algebraic structure given by the formula
$a_m+b_k\to a+b$; here $a_m\to a,b_k\to b$ for the
sequences $a_m\in M_m,b_k\in M_k$.  
The homomorphisms $\varphi_i$ can be arranged into  a graph as follows. 
Let  $M_i=M_{i_1}\oplus\dots\oplus M_{i_k}$ and 
$M_{i'}=M_{i_1'}\oplus\dots\oplus M_{i_k'}$ be 
the semi-simple $C^*$-algebras and $\varphi_i: M_i\to M_{i'}$ the  homomorphism. 
One has the two sets of vertices $V_{i_1},\dots, V_{i_k}$ and $V_{i_1'},\dots, V_{i_k'}$
joined by the $a_{rs}$ edges, whenever the summand $M_{i_r}$ contains $a_{rs}$
copies of the summand $M_{i_s'}$ under the embedding $\varphi_i$. 
As $i$ varies, one obtains an infinite graph called a {\it Bratteli diagram} of the
$AF$-algebra.

\medskip
{\bf E.}  
By ${\Bbb A}_{\theta}$ we denote an $AF$-algebra given by the Bratteli
diagram of Fig. 1. It is known that $K_0({\Bbb A}_{\theta})\cong {\Bbb Z}^2$
and $K_0^+({\Bbb A}_{\theta})=\{(p,q)\in {\Bbb Z}^2~|~p+\theta q\ge 0\}$. 
The $AF$-algberas ${\Bbb A}_{\theta}, {\Bbb A}_{\theta'}$ are stably
isomorphic, i.e. ${\Bbb A}_{\theta}\otimes {\cal K}\cong {\Bbb A}_{\theta'}\otimes {\cal K}$,  
if and only if ${\Bbb Z}+\theta {\Bbb Z}={\Bbb Z}+\theta'{\Bbb Z}$ as the subsets  of ${\Bbb R}$.

\subsection{The functor F}
{\bf A.} The Hubbard-Masur theory (\S 1.1) has been treated in a  general setting so far.
From now on, we  switch to the case $g=1$ (complex torus). Notice that $S=\widetilde S\cong T^2$,
since every  holomorphic quadratic differential $q$ on the complex torus is the square of
a holomorphic differential $\omega$.

\medskip
{\bf B.}
Let $\phi=Re~\omega$ be a 1-form defined by $\omega$. Since $\omega$
is holomorphic, $\phi$ is a closed $1$-form on $T^2$. The ${\Bbb R}$-isomorphism
$h_q: H^0(S,\Omega)\to Hom~(H_1(T^2); {\Bbb R})$, as explained,
is given by the formulas:
$$
\left\{
\begin{array}{cc}
\lambda_1 &= \int_{\gamma_1}\phi\\
\lambda_2 &= \int_{\gamma_2}\phi,
\end{array}
\right.
$$
where $\{\gamma_1,\gamma_2\}$ is a basis in the first homology group of $T^2$. 
We further assume that, after a proper choice of the basis,  $\lambda_1,\lambda_2$ are positive real
numbers.

\strut
\pagebreak

\medskip
{\bf C.}
Denote by $\Phi_{T^2}$ the space of measured foliations on $T^2$.
Each ${\cal F}\in \Phi_{T^2}$ is (measure) equivalent to a foliation
by a family of the parallel lines of a slope $\theta$ and the invariant (transverse)
measure $\mu$ (Fig.2).

\begin{figure}[here]
\begin{picture}(300,60)(-30,0)

\put(130,10){\line(1,0){40}}
\put(130,10){\line(0,1){40}}
\put(130,50){\line(1,0){40}}
\put(170,10){\line(0,1){40}}

\put(130,40){\line(2,1){20}}
\put(130,30){\line(2,1){40}}
\put(130,20){\line(2,1){40}}
\put(130,10){\line(2,1){40}}

\put(150,10){\line(2,1){20}}

\end{picture}

\caption{The measured foliation ${\cal F}$ on $T^2={\Bbb R}^2/{\Bbb Z}^2$.}
\end{figure}

\medskip\noindent
We use the notation ${\cal F}^{\mu}_{\theta}$ for such a foliation. 
There exists a simple  relationship between the reals  $(\lambda_1,\lambda_2)$
and $(\theta,\mu)$. Indeed, the closed $1$-form $\phi=Const$ defines a measured foliation,
${\cal F}^{\mu}_{\theta}$, so that
$$
\left\{
\begin{array}{ccc}
\lambda_1  &= \int_{\gamma_1}\phi &= \int_0^1\mu dx\\
\lambda_2  &= \int_{\gamma_2}\phi &= \int_0^1\mu dy
\end{array}
\right.
\hbox{,  where}
\quad
{dy\over dx}=\theta.
$$
By the integration:
$$
\left\{
\begin{array}{ccc}
\lambda_1  &= \int_0^1\mu dx &= \mu\\
\lambda_2  &= \int_0^1\mu\theta  dx &= \mu\theta.
\end{array}
\right.
$$
Thus, one gets $\mu=\lambda_1$ and $\theta={\lambda_2\over\lambda_1}$.

\medskip
{\bf D.}
Recall that the Hubbard-Masur theory (\S 1.1.C)  establishes a homeomorphism
$h: T_S(1)\to \Phi_{T^2}$, where $T_S(1)\cong {\Bbb H}=\{\tau: Im~\tau>0\}$
is the Teichm\"uller space of the  torus. Denote by $\omega_N$ an invariant
(N\'eron) differential of the complex torus ${\Bbb C}/(\omega_1{\Bbb Z}+\omega_2{\Bbb Z})$.
It is well known that $\omega_1=\int_{\gamma_1}\omega_N$ and 
$\omega_2=\int_{\gamma_2}\omega_N$, where $\gamma_1$ and $\gamma_2$ are the meridians of the torus.
Let $\pi$ be a projection acting by the formula $(\theta,\mu)\mapsto \theta$. 
An explicit formula for the functor  $F: {\cal E}\to {\cal A}$  is given by the composition:
$F=\pi\circ h$, where $h$ is the Hubbard-Masur homeomorphism. 
In other words, one gets the following (explicit) correspondence between the complex 
and noncommutative tori: 
$$
E_{\tau}=E_{(\int_{\gamma_2}\omega_N) / (\int_{\gamma_1}\omega_N)}
\buildrel\rm h\over
\longmapsto
{\cal F}^{\int_{\gamma_1}\phi}_{(\int_{\gamma_2}\phi)/(\int_{\gamma_1}\phi)}
\buildrel\rm\pi \over
\longmapsto
{\Bbb A}_{(\int_{\gamma_2}\phi)/(\int_{\gamma_1}\phi)}= {\Bbb A}_{\theta},
$$
where $E_{\tau}={\Bbb C}/({\Bbb Z}+\tau {\Bbb Z})$.

\section{Proof}
\subsection{Proof of lemma 1}
Let $\varphi: E_{\tau}\to E_{\tau'}$ be an isogeny of the elliptic curves. 
The action of $\varphi$ on the homology basis $\{\gamma_1,\gamma_2\}$
of $T^2$ is given by the formulas:
\begin{equation}\label{eq1}
\left\{
\begin{array}{cc}
\gamma_1' &= a\gamma_1+b\gamma_2\nonumber\\
\gamma_2' &= c\gamma_1+d\gamma_2
\end{array}
\right.
\hbox{,  where}
\left(\matrix{a & b\cr c & d}\small\right)\in M_2({\Bbb Z}). 
\end{equation}
Recall that the functor $F: {\cal E}\to {\cal A}$ is given by the formula:
\begin{equation}\label{eq2}
\tau={\int_{\gamma_2}\omega_N\over\int_{\gamma_1}\omega_N}
\longmapsto
\theta={\int_{\gamma_2}\phi\over\int_{\gamma_1}\phi},
\end{equation}
where $\omega_N$ is an invariant differential on $E_{\tau}$
and $\phi=Re~\omega$ is a closed 1-form on $T^2$.

\medskip
(i) From the left-hand side of (\ref{eq2}), one obtains 
\begin{equation}
\left\{
\begin{array}{ccccc}
\omega_1' &= \int_{\gamma_1'}\omega_N &=  \int_{a\gamma_1+b\gamma_2}\omega_N  &= 
a\int_{\gamma_1}\omega_N +b\int_{\gamma_2}\omega_N &= a\omega_1+b\omega_2\nonumber\\
\omega_2' &= \int_{\gamma_2'}\omega_N &=  \int_{c\gamma_1+d\gamma_2}\omega_N  &= 
c\int_{\gamma_1}\omega_N +d\int_{\gamma_2}\omega_N &= c\omega_1+d\omega_2,
\end{array}
\right.
\end{equation}
and therefore $\tau'={\int_{\gamma_2'}\omega_N\over\int_{\gamma_1'}\omega_N}=
{c+d\tau\over a+b\tau}$.

\medskip
(ii) From the right-hand side of (\ref{eq2}), one obtains
\begin{equation}
\left\{
\begin{array}{ccccc}
\lambda_1' &= \int_{\gamma_1'}\phi &=  \int_{a\gamma_1+b\gamma_2}\phi  &= 
a\int_{\gamma_1}\phi +b\int_{\gamma_2}\phi &= a\lambda_1+b\lambda_2\nonumber\\
\lambda_2' &= \int_{\gamma_2'}\phi &=  \int_{c\gamma_1+d\gamma_2}\phi  &= 
c\int_{\gamma_1}\phi +d\int_{\gamma_2}\phi &= c\lambda_1+d\lambda_2,
\end{array}
\right.
\end{equation}
and therefore $\theta'={\int_{\gamma_2'}\phi\over\int_{\gamma_1'}\phi}=
{c+d\theta\over a+b\theta}$.  Comparing (i) and (ii),  one gets the conclusion
of the first part of lemma \ref{lm1}. To prove the second part, recall that
the invertible isogeny is an isomorphism of the elliptic curves.
In this case  $\small\left(\matrix{a & b\cr c & d}\small\right)\in GL_2({\Bbb Z})$
and $\theta'=\theta ~mod~GL_2({\Bbb Z})$. Therefore $F$ sends the isomorphic
elliptic curves to the stably isomorphic noncommutative tori. The second part 
of lemma \ref{lm1} is proved. 

It follows from the proof that $F: {\cal E}\to {\cal A}$ is a covariant
functor. Indeed, $F$ preserves the morphisms and does not reverse the arrows:
$F(\varphi_1\varphi_2)=\varphi_1\varphi_2=F(\varphi_1)F(\varphi_2)$
for any pair of the isogenies $\varphi_1,\varphi_2\in Mor~({\cal E})$. 
$\square$

\subsection{Proof of Theorem 1}
The following lemma will be helpful.
\begin{lem}\label{lm2}
Let ${\goth m}\subset {\Bbb R}$ be a module of the rank 2, i.e
${\goth m}={\Bbb Z}\lambda_1+{\Bbb Z}\lambda_2$, where
$\theta={\lambda_2\over\lambda_1}\not\in {\Bbb Q}$. 
If ${\goth m}'\subseteq {\goth m}$ is a submodule of the rank 2,
then ${\goth m}'=k {\goth m}$, where either:

\medskip
(i) $k\in {\Bbb Z}-\{0\}$ and $\theta\in {\Bbb R}-{\Bbb Q}$, or 

\smallskip
(ii) $k$ and $\theta$ are the irrational numbers of  
a quadratic number field. 
\end{lem}
{\it Proof.} Any rank 2 submodule of $m$ can be written as  
${\goth m}'=\lambda_1'{\Bbb Z}+{\lambda_2'}{\Bbb Z}$, where
\begin{equation}
\left\{
\begin{array}{cc}
\lambda_1'  &= a\lambda_1 +b\lambda_2\\
\lambda_2'  &= c\lambda_1 +d\lambda_2
\end{array}
\right.
\qquad \hbox{and} \qquad 
\small\left(\matrix{a & b\cr c & d}\small\right)\in M_2({\Bbb Z}).
\end{equation}

\bigskip
(i) Let us assume that $b\ne 0$.  Let $\Delta= (a+d)^2-4(ad-bc)$ and $\Delta'=(a+d)^2-4bc$. 
We shall consider the following cases.

\medskip
{\bf Case 1: $\Delta>0$ and $\Delta\ne m^2$, $m\in {\Bbb Z}-\{0\}$.} 
The real number  $k$ can be determined  from the equations: 
\begin{equation}\label{eq8}
\left\{
\begin{array}{ccc}
\lambda_1'  &= k\lambda_1 &= a\lambda_1 +b\lambda_2\\
\lambda_2'  &= k\lambda_2 &= c\lambda_1 +d\lambda_2.
\end{array}
\right.
\end{equation}
Since  $\theta={\lambda_2\over\lambda_1}$, one gets the  equation $\theta={c\theta+d\over a\theta+b}$
by taking a ratio of the two  equations above. A quadratic equation for $\theta$ writes
as $b\theta^2+(a-d)\theta-c=0$. The discriminant of the equation coincides with $\Delta$
and therefore there exist the real roots $\theta_{1,2}={a-d\pm\sqrt{\Delta}\over 2c}$.
Moreover, $k=a+b\theta=a+{b\over 2c}(a-d\pm\sqrt{\Delta})$. Since $\Delta$ is
not the square of an integer, $k$ and $\theta$ are the irrationalities of the quadratic number field
${\Bbb Q}(\sqrt{\Delta})$.

\medskip
{\bf Case 2: $\Delta>0$ and $\Delta=m^2$, $m\in {\Bbb Z}-\{0\}$.}
Note that $\theta={a-d\pm |m|\over 2c}$ is a rational number. 
Since $\theta$ does not satisfy the rank assumption of the lemma,
the case should be omitted.

\medskip
{\bf Case 3: $\Delta=0$.}
The quadratic equation has a double root $\theta={a-d\over 2c}\in {\Bbb Q}$.
This case leads to a module of the rank $1$, which is contrary to an assumption of the
lemma.

\medskip
{\bf Case 4: $\Delta<0$ and $\Delta'\ne m^2$, $m\in {\Bbb Z}-\{0\}$.}
Let us  define  a new basis
$\{\lambda_1'',\lambda_2''\}$ in ${\goth m}'$ so that 
\begin{equation}
\left\{
\begin{array}{cc}
\lambda_1''  &=  \lambda_1'\\
\lambda_2''  &= -\lambda_2'.
\end{array}
\right.
\end{equation}
Then:
\begin{equation}
\left\{
\begin{array}{cc}
\lambda_1''  &= a\lambda_1 +b\lambda_2\\
\lambda_2''  &= -c\lambda_1 -d\lambda_2,
\end{array}
\right.
\end{equation}
and $\theta={\lambda_2''\over\lambda_1''}={-c-d\theta\over a+b\theta}$.
The quadratic equation for $\theta$ has the form $b\theta^2+(a+d)\theta+c=0$,
whose discrimimant is $\Delta'=(a+d)^2-4bc$.
Let us show that $\Delta'>0$. Indeed, $\Delta= (a+d)^2-4(ad-bc)<0$ 
and the evident inequality $-(a-d)^2\le 0$ have the same sign, and
we shall add them up. After an obvious elimination, one gets $bc<0$.
Therefore $\Delta'$ is a  sum of the two positive integers, which is   
always a positive integer. 
Thus,  there exist the real roots $\theta_{1,2}={-a-d\pm\sqrt{\Delta'}\over 2b}$.
Moreover, $k=a+b\theta={1\over 2}(a-d\pm\sqrt{\Delta'})$. Since $\Delta'$ is
not the square of an integer, $k$ and $\theta$ are the irrational numbers in the 
quadratic field ${\Bbb Q}(\sqrt{\Delta'})$.

\medskip
{\bf Case 5: $\Delta<0$ and $\Delta'=m^2$, $m\in {\Bbb Z}-\{0\}$.}
Note that $\theta={-a-d\pm |m|\over 2b}$ is a rational number.
Since $\theta$ does not satisfy the rank assumption of the lemma,
the case should be omitted.

\bigskip
(ii) Assume that $b=0$.

\medskip
{\bf Case 1: $a-d\ne 0$.}
The quadratic equation for $\theta$ degenerates to a 
linear eauation $(a-d)\theta+c=0$. The root $\theta={c\over d-a}\in {\Bbb Q}$
does not satisfy the rank assumption again, and we omit the case.

\medskip
{\bf Case 2: $a=d$ and $c\ne 0$.}
It is easy to see, that the set of the solutions for $\theta$ is an empty set.

\medskip
{\bf Case 3: $a=d$ and $c=0$.}
Finally, in this case all coefficients of the quadratic equation vanish,
so that any $\theta\in {\Bbb R}-{\Bbb Q}$ is a solution. Note that in the view of (\ref{eq8}),
$k=a=d\in {\Bbb Z}$. Thus, one gets case (i) of the lemma. 
Since there are no other possiblities left,  lemma \ref{lm2} is proved.
$\square$

\begin{lem}\label{lm3}
Let $E$ be an elliptic curve with a complex multiplication
and $h$ the Hubbard-Masur map, which acts by the formulas of \S 1.3.D. 
Consider a module  $h(Isom~(E))=\mu_E({\Bbb Z}+\theta_E{\Bbb Z}):= {\goth m}_E$.
Then:

\medskip
(i) $\theta_E$ is a quadratic irrationality, 

\smallskip
(ii) $\mu_E\in {\Bbb Q}$ (up to a choice of $h$).  
\end{lem}
{\it Proof.} (i) Since $E$ has a complex multiplication, $End~(E)> {\Bbb Z}$.
In particular, there exists a nontrivial endomorphism $\varphi$, i.e an endomorphism which is not 
the multiplication by $k\in {\Bbb Z}$.  By the lemma \ref{lm1}, $\varphi$ defines a submodule ${\goth m}_E'$
of the rank $2$ of the module  ${\goth m}_E$. By the lemma \ref{lm2}, ${\goth m}_E'=k {\goth m}_E$ for a
$k\in {\Bbb R}$. Since $\varphi$ is a nontrivial endomorphism, $k\not\in {\Bbb Z}$. 
Thus, the option (i) of lemma \ref{lm2} is excluded. Therefore,   
by the item (ii) of  lemma \ref{lm2},  $\theta_E$ is a quadratic irrationality.

\smallskip
(ii) Recall that $E_{\cal F}\subset Q- \{0\}$ is the space of holomorphic
differentials on the complex torus, whose  horizontal trajectory
structure is equivalent to given  measured foliation ${\cal F}={\cal F}^{\mu}_{\theta}$. 
We shall vary ${\cal F}_{\theta}^{\mu}$, thus varying the Hubbard-Masur homeomorphism
$h=h({\cal F}^{\mu}_{\theta}): E_{\cal F}\to T(1)$. Namely, consider a 1-parameter
continuous family of such maps  $h=h_{\mu}$, where $\theta=Const$ and $\mu\in {\Bbb R}$.  
Recall that $\mu_E=\lambda_1=\int_{\gamma_1}\phi$, where $\phi=Re~\omega$ and $\omega\in E_{\cal F}$. 
The family $h_{\mu}$ generates  a family $\omega_{\mu}=h^{-1}_{\mu}(C)$, where $C$ is a fixed
point in $T(1)$. Denote by $\phi_{\mu}$ and $\lambda_1^{\mu}$ the corresponding families
of the closed 1-forms and their periods, respectively. By the continuity, $\lambda_1^{\mu}$
takes on a rational value for a $\mu=\mu'$. (Actually, every  neighborhood of $\mu_0$
contains such a $\mu'$.) Thus, $\mu_E\in {\Bbb Q}$ for the Hubbard-Masur homeomorphism $h=h_{\mu'}$.  
$\square$

\bigskip
Lemma \ref{lm3} implies (i) of the theorem \ref{thm1}. To prove (ii),  notice that 
when  $E_1,E_2\in Isom (E_{CM})$, the 
respective modules ${\goth m}_1={\goth m}_2$. It follows from 
the fact that an isomorphism  between the elliptic curves  
corresponds to a change of basis in the module ${\goth m}$ (lemma \ref{lm1}).
Theorem \ref{thm1} is proved.
$\square$

\section{Arithmetic complexity of the noncommutative tori}
Let ${\Bbb A}_{\theta}$ be the noncommutative torus with a real multiplication.
Since $\theta$ is a quadratic irrationality,  the regular continued fraction of 
$\theta$ is eventually periodic:
\begin{equation}
\theta=[a_0, a_1,\dots, \overline{a_{k+1},\dots, a_{k+p}}],
\end{equation}
where $\overline{a_{k+1},\dots, a_{k+p}}$ is the minimal period of the
continued fraction. 
\begin{dfn}
Let us call the number $c({\Bbb A}_{\theta})=p$ an arithmetic complexity
of the noncommutative torus with real multiplication. 
\end{dfn}
\begin{lem}
The number $c({\Bbb A}_{\theta})$ is an invariant of the stable
isomorphism class of the noncommutative torus ${\Bbb A}_{\theta}$. 
\end{lem}
{\it Proof.} It follows from lemma \ref{lm1} that 
${\Bbb A}_{\theta}, {\Bbb A}_{\theta'}$  are stably isomorphic 
if and only if $\theta'=\theta~mod~SL(2,{\Bbb Z})$. By the 
main property of the regular continued fractions, the
expansion of $\theta$ and $\theta'$ must coincide, except possibly a finite
number of the entries. Since the continued fraction of $\theta$ is eventually
periodic, so must be the continued fraction of $\theta'$. Moreover,
the minimal periods of $\theta,\theta'$ must coincide as well as
their lengths. Thus $c({\Bbb A}_{\theta'})=c({\Bbb A}_{\theta})$.
$\square$ 
 
\begin{exm}
Let us find an  arithmetic complexity of the noncommutative
torus ${\Bbb A}_{3\sqrt{6}}$. The continued fraction expansion
of $3\sqrt{6}=\sqrt{54}$ is $[7; \overline{2,1,6,1,2,14}]$.
Since the continued fraction is six-periodic, 
we have $c({\Bbb A}_{3\sqrt{6}})=6$. 
\end{exm}
It is very useful to think of the normalized period 
$(1, {a_{k+2}\over a_{k+1}}, \dots, {a_{k+p}\over a_{k+1}})$
of ${\Bbb A}_{\theta}$ as coordinates of the `rational points'
of the noncommutative torus, taken up to a cyclic permutation. 
In the sense, such points are the generators of an abelian group
of all rational points of ${\Bbb A}_{\theta}$ modulo the points
of a finite order.  
\begin{cnj}
$c({\Bbb A}_{\theta_{E_{CM}}})=rk~(E_{CM})+1$. 
\end{cnj}

\bigskip\noindent
{\sf Acknowledgments.} 
I am grateful to Lawrence ~D. ~Taylor for the helpful remarks regarding  the 
manuscript. The referee's suggestions are kindly acknowledged and 
incorporated in the text.



\vskip1cm

\textsc{The Fields Institute for Mathematical Sciences, Toronto, ON, Canada,  
E-mail:} {\sf igor.v.nikolaev@gmail.com}

\end{document}